\input ./style/arxiv-general.cfg
\documentclass[aos,MSNbibl,nameyear,seceqn,dvips]{arximspdf}
\makeatletter
   \@ifpackageloaded{graphicx}{}{\usepackage{graphicx}}
\makeatother
\usepackage{mathrsfs}

\doi{10.1214/15-AOS1376}
\volume{44}
\issue{2}
\pubyear{2016}
\firstpage{540}
\lastpage{563}
\docsubty{FLA}

\makeatletter
\newcommand{\rrVert}{\Vert}
\newcommand{\llVert}{\Vert}
\newcommand{\overset}{\stackrel}
\newproclaim{definition}{Definition}
\newtheorem{theorem}{Theorem}

\renewcommand{\mid}{|}

\newcommand{\off}{\operatorname{off}}
\newcommand{\proj}{\mathrm{Proj}}
\newcommand{\vc}{\operatorname{vec}}
\newcommand{\tr}{\operatorname{trace}}
\newcommand{\dd}{\mathrm{d}}
\newcommand{\eye}{\mathbf{I}}
\makeatother

\begin{document}
\begin{frontmatter}

\title{Estimating multivariate latent-structure models}
\runtitle{Multivariate latent-structure models}

\begin{aug}
\author[A]{\fnms{St\'ephane}~\snm{Bonhomme}\thanksref{M1,T1}\ead[label=e1]{sbonhomme@uchicago.edu}},
\author[B]{\fnms{Koen}~\snm{Jochmans}\corref{}\thanksref{M2,T2}\ead[label=e2]{koen.jochmans@sciencespo.fr}}
\and
\author[C]{\fnms{Jean-Marc}~\snm{Robin}\thanksref{M2,M3,T3}\ead[label=e3]{jeanmarc.robin@sciencespo.fr}}
\runauthor{S. Bonhomme, K. Jochmans and J.-M. Robin}
\affiliation{University of Chicago\thanksmark{M1}, Sciences
Po\thanksmark{M2} and University College London\thanksmark{M3}}
\address[A]{S. Bonhomme\\
Department of Economics\\
University of Chicago\\
1126 E.~59th Street\\
Chicago, Illinois 60637\\
USA\\
\printead{e1}}
\address[B]{K. Jochmans\\
Department of Economics\\
Sciences Po\\
28 rue des Saints P\`eres\\
75007 Paris\\
France\\
\printead{e2}}
\address[C]{J.-M. Robin\\
Department of Economics\\
Sciences Po\\
28 rue des Saints P\`eres\\
75007 Paris\\
France\\
and\\
Department of Economics\\
University College London\\
Drayton House\\
30 Gordon Street\\
London WC1 H0AX\\
United Kingdom\\
\printead{e3}}
\end{aug}
\thankstext{T1}{Supported by European Research Council Grant
ERC-2010-StG-0263107-ENMUH.}
\thankstext{T2}{Supported by Sciences Po's SAB grant ``Nonparametric
estimation of finite mixtures.''}
\thankstext{T3}{Supported by European Research Council Grant
ERC-2010-AdG-269693-WASP and by Economic and Social Research Council
Grant RES-589-28-0001 through the Centre for Microdata Methods and Practice.}

%
\received{\smonth{5} \syear{2015}}
%
\revised{\smonth{8} \syear{2015}}

%
\begin{abstract}
A constructive proof of identification of multilinear decompositions of
multiway arrays is presented. It can be applied to show identification
in a variety of multivariate latent structures. Examples are
finite-mixture models and hidden Markov models. The key step to show
identification is the joint diagonalization of a set of matrices in the
same nonorthogonal basis. An estimator of the latent-structure model
may then be based on a sample version of this joint-diagonalization
problem. Algorithms are available for computation and we derive
distribution theory. We further develop asymptotic theory for
orthogonal-series estimators of component densities in mixture models
and emission densities in hidden Markov models.
\end{abstract}

%
\end{frontmatter}

\setcounter{footnote}{3}

\section{Introduction}
Latent structures are a popular tool for modeling the dependency
structure in multivariate data. Two important examples are
finite-mixture models [see \citet{McLachlanPeel2000}] and hidden
Markov models [see \citet{CappeMoulinesRyden2005}]. Although
these models arise frequently in applied work, the question of their
nonparametric identifiability has attracted substantial attention only
quite recently. \citet{AllmanMatiasRhodes2009} used algebraic results
on the uniqueness of decompositions of multiway arrays due to \citeauthor{Kruskal1976}
(\citeyear{Kruskal1976,Kruskal1977}) to establish identification in a variety of
multivariate latent-structure models. Their setup covers both finite
mixtures and hidden Markov models, among other models, and their
findings substantially generalize the earlier work of \citet{Green1951,Anderson1954,Petrie1969}, \citet
{HettmanspergerThomas2000}, \citet{HallZhou2003}, and \citet
{HallNeemanPakyariElmore2005}.

Despite these positive identification results, direct application of
Kruskal's method does not provide an estimator. Taking identification
as given, some authors have developed EM-type approaches to
nonparametrically estimate both multivariate finite mixtures
[\citet{BenagliaChauveauHunter2009}; Levine, Hunter and Chauveau (\citeyear{LevineHunterChauveau2011})] and hidden Markov models [\citet
{GassiatCleynenRobin2013}]. Numerical studies suggest that these
estimators are well behave. However, their statistical properties---their consistency, convergence rates, and asymptotic
distribution---are difficult to establish and are currently unknown.\footnote{There
are results on inference in semi- and nonparametric finite-mixture
models and hidden Markov models in several more restrictive settings.
These include location models [\citet
{BordesMotteletVandekerkhove2006}; \citet
{HunterWangHettmansperger2007}; and \citet{GassiatRousseau2014}],
multivariate finite mixtures with identically distributed outcome
variables [\citet{HettmanspergerThomas2000}; \citet
{BonhommeJochmansRobin2014}], and two-component mixtures [\citet
{HallZhou2003}; \citet{HenryJochmansSalanie2013}].}


In this paper, we show that the multilinear structure underlying the
results of \citet{AllmanMatiasRhodes2009} can be used to obtain a
constructive proof of identification in a broad class of
latent-structure models. We show that the problem of decomposing a
multiway array can be reformulated as the problem of simultaneously
diagonalizing a collection of matrices. This is a least-squares problem
that has received considerable attention in the literature on
independent component analysis and blind source separation [see
\citet{ComonJutten2010}]. Moreover, algorithms exist to recover
the joint diagonalizer in a
computationally efficient
manner; see \citet
{FuGao2006}, \citeauthor{IferroudjeneAbedMeraimBelouchrani2009} (\citeyear{IferroudjeneAbedMeraimBelouchrani2009,IferroudjeneAbed-MeraimBelouchrani2010})
and \citeauthor{LucianiAlbera2010} (\citeyear{LucianiAlbera2010,LucianiAlbera2014}).

We propose estimating the parameters of the latent-structure model by
solving a sample version of the simultaneous diagonalization problem.
We provide distribution theory for this estimator below. Under weak
conditions, it converges at the parametric rate and is asymptotically
normal. Using this result, we obtain estimators of finite-mixture
models and hidden Markov models that have standard asymptotic
properties. Moreover, the fact that the dependency structure in the
data is latent does not translate into a decrease in the convergence
rate of the estimators. As such, this paper is the first to derive the
asymptotic behavior of nonparametric estimators of multivariate
finite-mixture models of the form defined in \citet{HallZhou2003} for
more than two latent classes and of hidden Markov models of the form in
\citet{GassiatCleynenRobin2013}. Furthermore, our approach can be
useful in the analysis of random graph models [\citet
{AllmanMatiasRhodes2011}] and stochastic blockmodels [\citet
{SnijdersNowicki1997}; \citet{RoheChatterjeeYu2011}], although we
do not consider such models in detail in this paper. In a simulation
study, we find that our approach performs well in small samples.

There is a large literature on parallel factor analysis and canonical
polyadic decompositions of tensors building on the work of \citeauthor{Kruskal1976} (\citeyear{Kruskal1976,Kruskal1977}); see, for example,
\citet{DeLathauwerDeMoorVandewalle2004},
\citet{DeLathauwer2006},
\citeauthor{DomanovDeLathauwer2013a} (\citeyear{DomanovDeLathauwer2013a,DomanovDeLathauwer2013b,DomanovDeLathauwer2014b,DomanovDeLathauwer2014a}),
\citet{AnandKumarGeHsuKakadeTelgarsky2014} and
\citeauthor{ChiantiniOttavianiVannieuwenhoven2014} (\citeyear{ChiantiniOttavianiVannieuwenhoven2014,ChiantiniOttavianiVannieuwenhoven2015}).
Although our strategy has some similarity with this literature, both
our conclusions and our simultaneous diagonalization problem are
different. Most importantly, our simultaneous diagonalization
formulation can deal with noise, making it useful as a tool for
statistical inference.

In the context of multivariate finite mixtures of identically
distributed variables, \citet{KasaharaShimotsu2009} and \citet
{BonhommeJochmansRobin2014} also used (different) joint-diagonalization
arguments to obtain nonparametric identification results. However, the
approaches taken there are different from the one developed in this
paper and cannot be applied as generally.

We start out by motivating our approach via a discussion on the
algebraic structure of multivariate finite-mixture models and hidden
Markov models. We then present our identification strategy in a generic
setting. After this we turn to estimation and inference, and to the
development of asymptotic theory. Next, the theory is used to set up
orthogonal-series estimators of component densities in a finite-mixture
model, and to show that these have the standard univariate convergence
rates of series estimators. Finally, the orthogonal-series density
estimator is put to work in simulation experiments involving finite
mixtures and a hidden Markov model. The supplementary material
[\citet{BonhommeJochmansRobin2014b}] contains some additional
results and discussion, as well as all technical proofs.

\section{Motivating examples} \label{secexamples}
We start by introducing three examples to motivate our subsequent developments.

\subsection{Finite-mixture models for discrete measurements}\label{sec2.1}
Let $Y_1,Y_2,\ldots,Y_q$ be observable random variables that are
assumed independent conditional on realizations of a latent random
variable $Z$. Suppose that $Z$ has a finite state space of known
cardinality $r$, which we set to $\lbrace1,2,\ldots,r \rbrace$
without loss of generality. Let $\bolds{\pi}=(\pi_1,\pi_2,\ldots
,\pi_r)^\prime$ be the probability distribution of $Z$, so $\pi_j>0$
and $\sum_{j=1}^r \pi_j=1$. Then the\vspace*{1pt} probability distribution of
$Y_1,Y_2,\ldots,Y_q$ is a multivariate finite mixture with mixing
proportions $\pi_1,\pi_2,\ldots,\pi_r$. The parameters of interest
are the mixing proportions and the distributions of $Y_1,Y_2,\ldots
,Y_q$ given $Z$. The $Y_i$ need not be identically distributed, so the
model involves $qr$ such conditional distributions.

Suppose that the scalar random variable $Y_i$ can take on a finite
number $\kappa_i$ of values. Let $\mathbf
{p}_{ij}=(p_{ij1},p_{ij2},\ldots,p_{ij\kappa_i})^\prime$ denote the
probability distribution of $Y_i$ given $Z=j$. Let $\bigotimes$ denote
the outer (tensor) product. The joint probability distribution of
$Y_1,Y_2,\ldots,Y_q$ given $Z=j$ then is the $q$-way table
\[
\bigotimes_{i=1}^q \mathbf{p}_{ij}=
\mathbf{p}_{1j} \otimes\mathbf {p}_{2j} \otimes\cdots\otimes
\mathbf{p}_{qj},
\]
which is of dimension $\kappa_1\times\kappa_2\times\cdots\times
\kappa_q$. The outer-product representation follows from the
conditional-independence restriction. Hence, the marginal probability
distribution of $Y_1,Y_2,\ldots,Y_q$ equals
%
\begin{equation}
\label{eqdiscretemixture} \mathbb{P} = \sum_{j=1}^r
\pi_j \bigotimes_{i=1}^q
\mathbf{p}_{ij},
\end{equation}
which is an $r$-linear decomposition of a $q$-way array. The parameters
of the mixture model are all the vectors making up the outer-product
arrays, $\lbrace\mathbf{p}_{ij}\rbrace$ and the coefficients of the
linear combination, $\lbrace\pi_j\rbrace$, transforming the
conditional distributions into the marginal distribution $\mathbb{P}$.

The $r$-linear decomposition is not restricted to the contingency
table. Indeed, any linear functional of $\mathbb{P}$ admits a
decomposition in terms of the same functional of the $\mathbf
{p}_{ij}$. Moreover, for any collection of vector-valued
transformations $y\mapsto\bolds{\chi}_i(y)$ we have
%
\begin{equation}
\label{eqlinfunc} E \Biggl[\bigotimes_{i=1}^q
\bolds{\chi}_i(Y_i) \Biggr] = \sum
_{j=1}^r \pi_j \bigotimes
_{i=1}^q E\bigl[\bolds{\chi }_i(Y_i)
\mid Z=j\bigr],
\end{equation}
provided the expectation exists. Of course, identification of linear
functionals follows from identification of the component distributions,
but (\ref{eqlinfunc}) can be useful for the construction of
estimators. To illustrate this, we turn to a model with continuous outcomes.

\subsection{Finite-mixture models for continuous measurements}
Suppose now that the $Y_i$ are continuously distributed random
variables. Let $f_{ij}$ be the density of $Y_i$ given $Z=j$. In this
case, the $q$-variate finite-mixture model with $r$ latent classes
states that the joint density function of the outcomes $Y_1,Y_2,\ldots
,Y_q$ factors as
%
\begin{equation}
\label{eqcontinuousmixture} \sum_{j=1}^r
\pi_j \prod_{i=1}^q
f_{ij},
\end{equation}
again for mixing proportions $\pi_1,\pi_2,\ldots,\pi_r$. This is an
infinite-dimensional version of (\ref{eqdiscretemixture}). Setting
$\bolds{\chi}_i$ in (\ref{eqlinfunc}) to a set of indicators that
partition the state space of $Y_i$ yields a decomposition as in (\ref
{eqdiscretemixture}) for a discretized version of the mixture model.
This approach has been used by \citet{AllmanMatiasRhodes2009} and \citet
{KasaharaShimotsu2014} in proving identification.

An alternative approach, which will prove convenient for the
construction of density estimators, is as follows. Suppose that
$(Y_1,Y_2,\ldots,Y_q)$ lives in the \mbox{$q$-}dimensional space $\mathscr
{Y}^q\subseteq\mathscr{R}^q$. Let $L_\rho^2[\mathscr{Y}]$ be the
space of functions that are square-integrable with respect to the
weight function $\rho$ on $\mathscr{Y}$, endowed with the inner product
\[
\langle h_1,h_2 \rangle= \int_{\mathscr{Y}}
h_1(y) h_2(y) \rho (y) \,\dd y,
\]
and the $L_\rho^2$-norm $\llVert  h \rrVert _2 = \sqrt
{\langle h,h \rangle}$. Let $\lbrace\varphi_k, k>0 \rbrace$ be a
class of functions that form a complete orthonormal basis for $L_\rho
^2[\mathscr{Y}]$. When $\mathscr{Y}$ is compact, polynomials such as
those belonging to the Jacobi class (e.g., Chebychev or Legendre polynomials) can serve this purpose. When $\mathscr{Y}=(-\infty
,+\infty)$, Hermite polynomials are a natural choice.

Assume that $f_{ij}\in L_\rho^2[\mathscr{Y}]$. The projection of
$f_{ij}$ onto the subspace spanned by $\varphi_1,\varphi_2,\ldots
,\varphi_{\varkappa}$ for any integer $\varkappa$ is
\[
\proj_{\varkappa} f_{ij} = \sum_{k=1}^\varkappa
b_{ijk} \varphi_k,
\]
where the
\[
b_{ijk} = \langle\varphi_k, f_{ij} \rangle= E
\bigl[{\varphi }_{k}(Y_i)\rho(Y_i)\mid Z=j
\bigr]
\]
are the (generalized) Fourier coefficients of $f_{ij}$. The projection
converges to $f_{ij}$ in $L^2_\rho$-norm, that is,
$
\llVert \proj_{\varkappa} f_{ij}-f_{ij} \rrVert _2
\rightarrow0
$
as $\varkappa\rightarrow\infty$. Such projections are commonly-used
tools in the approximation of functions and underlie orthogonal-series
estimators of densities.

The Fourier coefficients are not directly observable. For chosen
integers $\kappa_1,\kappa_2,\ldots,\kappa_q$, define
\[
\mathbf{b}_{ij} = E\bigl[\bolds{\varphi}_{\kappa_i}(Y_i)
\rho(Y_i)\mid Z=j\bigr],
\]
where
$
\bolds{\varphi}_{\kappa_i} =
(\varphi_1,\varphi_2,\ldots,\varphi_{\kappa_i})^\prime
$,
which are linear functionals of the $f_{ij}$. Then (\ref{eqlinfunc}) yields
%
\begin{equation}
\label{eqfouriersystem} \mathbb{B} = \sum_{j=1}^r
\pi_j \bigotimes_{i=1}^q
\mathbf{b}_{ij}
\end{equation}
for $\mathbb{B}=E[\bigotimes_{i=1}^q \bolds{\varphi}_{\kappa
_i}(Y_i)\rho(Y_i)]$. The latter expectation is a $q$-way array that
can be computed directly from the data. It contains the leading Fourier
coefficients of the $q$-variate density function of the data. Again,
the array $\mathbb{B}$ factors into a linear combination of multiway
arrays. In Section~\ref{secmixtures}, we will use this representation
to derive orthogonal-series density estimators that have standard
large-sample properties.

\subsection{Hidden Markov models} \label{subsechmm}
Let $\lbrace Y_i,Z_i \rbrace_{i=1}^q$ be a stationary sequence. $Z_i$
is a latent variable with finite state space $\lbrace1,2,\ldots,r
\rbrace$, for known $r$, and has first-order Markov dependence. Let
$\bolds{\pi}=(\pi_1,\pi_2,\ldots,\pi_r)^\prime$ be the
stationary distribution of $Z_i$. Write $\mathbf{K}$ for the $r\times
r$ matrix of transition probabilities; so $\mathbf{K}(j_1,j_2)$ is the
probability of moving from state $j_1$ to state $j_2$. The observable
scalar random variables $Y_1,Y_2,\ldots,Y_q$ are independent
conditional on realizations of $Z_1,Z_2,\ldots,Z_q$, and the
distribution of $Y_i$ only depends on the realization of $Z_i$. This is
a hidden Markov model with $r$ latent states and $q$ observable outcomes.

Suppose that $Y_i$ is discrete and that its state space contains
$\kappa$ points of support. Write $\mathbf{p}_j$ for the probability
vector of $Y_i$ given $Z_i=j$, that is, the emission distributions. Let
$\mathbf{P}=(\mathbf{p}_1,\mathbf{p}_2,\ldots,\mathbf{p}_r)$ be
the $\kappa\times r$ matrix of emission distributions and write
$\bolds{\Pi}=\operatorname{diag}(\pi_1,\pi_2,\ldots,\pi_r)$. The
Markovian assumption implies that $Y_i$ and $Z_{i-1}$ are independent
given $Z_i$. Hence, the columns of the matrix
\[
\mathbf{B} = \mathbf{P}\mathbf{K}^\prime= (\mathbf{b}_1,
\mathbf {b}_2,\ldots,\mathbf{b}_r)
\]
contain the probability distributions of $Y_i$ for given values of
$Z_{i-1}$. Likewise, $Y_i$ and $Z_{i+1}$ are independent given $Z_i$,
and so the matrix
\[
\mathbf{A} = \mathbf{P}\bolds{\Pi} \mathbf{K} \bolds{\Pi}^{-1} = (
\mathbf{a}_1,\mathbf{a}_2,\ldots,\mathbf{a}_r)
\]
gives the distributions of $Y_i$ for given values of $Z_{i+1}$.
Finally, $Y_{i-1}$, $Y_i$, and $Y_{i+1}$ are independent given $Z_i$.
Thus, with $q=3$ measurements, the hidden Markov model implies that the
contingency table of $(Y_1,Y_2,Y_3)$ factors as
%
\begin{equation}
\label{eqhmm} \mathbb{P} = \sum_{j=1}^r
\pi_j (\mathbf{a}_j\otimes\mathbf {p}_j
\otimes\mathbf{b}_j).
\end{equation}
A detailed derivation is provided in the supplementary material [\citet{BonhommeJochmansRobin2014b}]; also
see [\citet{GassiatCleynenRobin2013}, Theorem~2.1]
and [\citet{AllmanMatiasRhodes2009}, Section~6.1] for alternative
derivations. When $q>3$, we may bin several outcomes together and
proceed as before, by using the unfolding argument in Section~\ref{subsecunfolding}.

Equation (\ref{eqhmm}) shows that appropriate conditioning allows
viewing the hidden Markov model as a finite-mixture model, thus casting
it into the framework of finite mixtures with conditionally-independent
(although not identically-distributed) outcomes as in (\ref
{eqdiscretemixture}). Here, the parameters of interest are the
emission distributions $\lbrace\mathbf{p}_j \rbrace_{j=1}^r$ and the
stationary distribution of the Markov chain $\bolds{\pi}$, and also
the matrix of transition probabilities $\mathbf{K}$.

When the $Y_i$ are continuously distributed, (\ref{eqhmm}) becomes a
mixture as in (\ref{eqcontinuousmixture}), and we may again work with
projections of the densities onto an orthogonal basis.

\section{Algebraic structure and identification} \label{secid}
Our approach can be applied to $q$-variate structures that decompose as
$q$-ads, which are defined as follows.

\begin{definition} \label{defqad}
A $q$-dimensional array $\mathbb{X}\in\mathscr{R}^{\kappa_1\times
\kappa_2\times\cdots\times\kappa_q}$ is a $q$-ad if it can be
decomposed as
%
\begin{equation}
\label{eqqad} \mathbb{X} = \sum_{j=1}^r
\pi_j \bigotimes_{i=1}^q
\mathbf{x}_{ij}
\end{equation}
for some integer $r$, nonzero weights $\pi_1,\pi_2,\ldots,\pi_r$,
and vectors $\mathbf{x}_{ij}\in\mathscr{R}^{\kappa_i\times1}$.
\end{definition}

Our interest lies in nonparametrically recovering $\lbrace\mathbf
{x}_{ij}\rbrace$ and $\lbrace\pi_j\rbrace$ from knowledge of
$\mathbb{X}$ and $r$. Clearly, these parameters are not unique, in
general. For example, a permutation of the $\mathbf{x}_{ij}$ and $\pi
_j$ leaves $\mathbb{X}$ unaffected, and a common scaling of the
$\mathbf{x}_{ij}$ combined with an inverse scaling of the $\pi_j$,
too, does not change the $q$-way array. However, the work of
\citeauthor{Kruskal1976} (\citeyear{Kruskal1976,Kruskal1977}), \citet{SidiropoulosBro2000}, \citet
{JiangSidiropoulos2004} and \citeauthor{DomanovDeLathauwer2013a} (\citeyear{DomanovDeLathauwer2013a,DomanovDeLathauwer2013b}), among others, gives
simple sufficient conditions for uniqueness of the decomposition up to
these two indeterminacies. These conditions cannot be satisfied when $q<3$.

While permutational equivalence of possible decompositions of $\mathbb
{X}$ is an inherently unresolvable ambiguity, indeterminacy of the
scale of the vectors $\mathbf{x}_{ij}$ is undesirable in many
situations. Indeed, in arrays of the general form in (\ref
{eqlinfunc}), recovering the scale of the $\mathbf{x}_{ij}$ and the
constants $\pi_j$ is fundamental. In some cases, natural scale
restrictions may be present. Indeed, in (\ref{eqdiscretemixture}) the
$\mathbf{x}_{ij}$ are known to be probability distributions, and so
they have nonnegative entries that sum to one. Suitably combining
these restrictions with Kruskal's theorem, \citet
{AllmanMatiasRhodes2009} derived conditions under which the parameters
in finite mixtures and hidden Markov models are uniquely determined up
to relabelling of the latent classes.

We follow a different route to determine $q$-adic decompositions up to
permutational equivalence that does not require knowledge of the scale
of the $\mathbf{x}_{ij}$. We require that, apart from the $q$-way
array $\mathbb{X}$, lower-dimensional submodels are also observable.
By lower-dimensional submodels we mean arrays that factor as
%
\begin{equation}
\label{eqsubmodel} \sum_{j=1}^r
\pi_j \bigotimes_{i\in\mathcal{Q}}
\mathbf{x}_{ij}
\end{equation}
for sets $\mathcal{Q}$ that are subsets of the index set $\lbrace
1,2,\ldots,q\rbrace$. This is not a strong requirement in the models
we have in mind. For example, in the mixture model in~(\ref
{eqdiscretemixture}), lower-dimensional submodels are just the
contingency tables of a subset of the outcome variables. There, going
from a $q$-way table down to a $(q-1)$-table featuring all but the
$i$th outcome boils down to summing the array in the $i$th direction.
In more general situations, such as (\ref{eqlinfunc}) and in the
multilinear equation involving Fourier coefficients in particular, the
advantage of working with submodels over marginalizations of the model
is apparent. Indeed, in contrast to when the array is a contingency
table, here, there is no natural scale constraint on the $\mathbf
{x}_{ij}$. So, summing the array in one direction does not yield an
array that decomposes as in (\ref{eqsubmodel}). Nonetheless,
expectations concerning any subset of the random variables can still be
computed in (\ref{eqlinfunc}) and so submodels as defined in (\ref
{eqsubmodel}) are observable. In the supplementary material [\citet{BonhommeJochmansRobin2014b}] we adapt
our main identification result (Theorem~\ref{thmid1} below) to
settings where submodels are not available and marginalizations are
used instead.

Note that, throughout, we take $r$ in (\ref{eqqad}) to be known. This
ensures $\lbrace\mathbf{x}_{ij}\rbrace$ and $\lbrace\pi_j \rbrace
$ to be unambiguously defined. For a different $r$, there may exist a
different set of weights and vectors so that $\mathbb{X}$ factors as a
$q$-ad. The rank of $\mathbb{X}$ is the smallest integer $r$ needed to
arrive at a decomposition as in Definition \ref{defqad}. For example,
in the multivariate mixture model in Section~\ref{sec2.1}, $r$ is the number of
fitted mixture components and the rank is the smallest number of
components that would allow us to write the joint distribution of the
variables as a mixture that satisfies the required
conditional-independence restriction as in (\ref{eqdiscretemixture}).
The rank need not be equal to $r$. Moreover, besides the factorization
of $\mathbb{P}$ in terms of $\pi_1,\pi_2,\ldots,\pi_r$ and
$\lbrace\mathbf{p}_{i1},\mathbf{p}_{i2},\ldots,\mathbf{p}_{ir}
\rbrace$ in (\ref{eqdiscretemixture}), there may exist a different
set of, say, $r^\prime$ weights $\pi_1^\prime,\pi_2^\prime,\ldots
,\pi^\prime_{r^\prime}$ and distributions $\lbrace\mathbf
{p}_{i1}^\prime,\mathbf{p}_{i2}^\prime,\ldots,\mathbf
{p}_{ir^\prime}^\prime\rbrace$ that also yield a representation of
$\mathbb{P}$ as a mixture. Identifying the number of components is a
difficult issue. Recent work by \citet{KasaharaShimotsu2014} shows that
a simple lower bound on the number of components is nonparametrically
identified (and estimable).

\subsection{Unfolding} \label{subsecunfolding}
We can state our main identification result for three-way arrays
without loss of generality. This is so because any $q$-way array can be
unfolded into a $(q-1)$-way array, much like any matrix can be
transformed into a vector using the vec operator. Indeed, in any
direction $i\in\lbrace1,2,\ldots,q\rbrace$, a \mbox{$q$-}way array of
dimension $\kappa_1\times\kappa_2\times\cdots\times\kappa_q$ is
a collection of $\kappa_i$ $(q-1)$-way arrays, each of dimension
$\kappa_1\times\kappa_2\times\cdots\times\kappa_{i-1}\times\kappa
_{i+1}\times\cdots\times\kappa_q$. This collection can be stacked in any
of $i^\prime\in\lbrace1,2,\ldots,i-1,i+1,\ldots,q\rbrace$
directions, that is, $(q-1)$ different ways, to yield a $(q-1)$-way
array whose dimension will be $\kappa_1\times\kappa_2 \times\kappa
_{i}\kappa_{i^\prime}\times\cdots\times \kappa_q$. This unfolding process
can be iterated until it yields a three-way array. To write this
compactly, let $\bigodot$ be the Khatri--Rao product. Then, for vectors
$\mathbf{a}_1, \mathbf{a}_2, \ldots,\mathbf{a}_q$,
\[
\bigodot_{i=1}^q \mathbf{a}_i =
\mathbf{a}_1 \odot\mathbf{a}_2 \odot\cdots\odot
\mathbf{a}_q
\]
is the vector containing all interactions between the elements of the
$\mathbf{a}_i$. The end result of iterated unfolding toward direction
$i$, say, is a three-way array of the form
\[
\sum_{j=1}^r \pi_j \Bigl(
\bigodot_{i_1\in\mathcal{Q}_1} \mathbf{x}_{i_1j} \otimes \mathbf{x}_{ij}
\otimes \bigodot_{i_2\in\mathcal{Q}_2} \mathbf{x}_{i_2j} \Bigr),
\]
where $\mathcal{Q}_1$ and $\mathcal{Q}_2$ are two index sets that
partition $\lbrace1,2,\ldots,q \rbrace\setminus\lbrace i \rbrace
$. We will illustrate this in the context of density estimation in
Section~\ref{secmixtures}.

\subsection{Identification via simultaneous diagonalization}
We thus focus on a three-way array $\mathbb{X}$ of dimension $\kappa
_1\times\kappa_2\times\kappa_3$ that factors as a tri-ad, that is,
\[
\mathbb{X} = \sum_{j=1}^r
\pi_j (\mathbf{x}_{1j}\otimes\mathbf {x}_{2j}
\otimes\mathbf{x}_{3j}).
\]
Let $\mathbf{X}_i=(\mathbf{x}_{i1},\mathbf{x}_{i2},\ldots,\mathbf
{x}_{ir})$ and $\bolds{\Pi}=\operatorname{diag}(\pi_1,\pi_2,\ldots,\pi
_r)$. Also, for each pair $(i_1,i_2)$ with $i_1<i_2$ in $\lbrace1,2,3
\rbrace^2$, let
\[
\mathbb{X}_{\lbrace i_1,i_2 \rbrace} = \sum_{j=1}^r
\pi_j (\mathbf{x}_{i_1j}\otimes\mathbf{x}_{i_2j}).
\]
Note that, from (\ref{eqsubmodel}), $\mathbb{X}_{\lbrace i_1,i_2
\rbrace}$ is the lower-dimension submodel obtained from $\mathbb{X}$
by omitting the index $i_3$.

Our first theorem concerns identification of the $\mathbf{X}_i$ as the
eigenvalues of a set of matrices and is the cornerstone of our
argument. The proof of this result is constructive and will be the
basis for our estimator in Section~\ref{secestimation} below.

\begin{theorem}[(Columns of $\mathbf{X}_i$)] \label{thmid1}
If $\mathbf{X}_{i_1}$ and $\mathbf{X}_{i_2}$ both have full column
rank and $\mathbb{X}_{\lbrace i_1,i_2 \rbrace}$ is observable, then
$\mathbf{X}_{i_3}$ is identified up to a permutation matrix if all its
columns are different.
\end{theorem}

\begin{pf}
Without loss of generality, f\mbox{i}x $(i_1,i_2,i_3)=(1,2,3)$
throughout the proof. In each direction $i$, the three-way array
$\mathbb{X}$ consists of a collection of $\kappa_i$ matrices. Let
$\mathbf{A}_1,\mathbf{A}_2,\ldots,\mathbf{A}_{\kappa_3}$ denote
these matrices for $i=3$. So, the matrix $\mathbf{A}_k$ is obtained
from $\mathbb{X}$ by fixing its third index to the value $k$, that is,
$\mathbf{A}_k = \mathbb{X}(:,:,k)$, using obvious array-indexing
notation. Also, let $\mathbf{A}_0 = \mathbb{X}_{\lbrace1,2\rbrace
}$. Note that all of $\mathbf{A}_0$ and $\mathbf{A}_1,\mathbf
{A}_2,\ldots,\mathbf{A}_{\kappa_3}$ are observable matrices of
dimension $\kappa_1\times\kappa_2$.

The lower-dimensional submodel $\mathbf{A}_0$ has the structure
%
\begin{equation}
\label{eqsub} \mathbf{A}_0 = \mathbf{X}_1 \bolds{
\Pi} \mathbf{X}_2^\prime.
\end{equation}
Because the matrices $\mathbf{X}_1$ and $\mathbf{X}_2$ both have rank
$r$ and because all $\pi_j$ are nonzero by definition, the matrix
$\mathbf{A}_0$, too, has rank $r$. Therefore, it has a singular-value
decomposition
\[
\mathbf{A}_0 = \mathbf{U}\mathbf{S}\mathbf{V}^\prime
\]
for unitary matrices $\mathbf{U}$ and $\mathbf{V}$ of dimension
$\kappa_1\times r$ and $\kappa_2\times r$, respectively, and a
nonsingular $r\times r$ diagonal matrix $\mathbf{S}$. Now construct
$\mathbf{W}_1=\mathbf{S}^{-1/2}\mathbf{U}^\prime$ and $\mathbf
{W}_2=\mathbf{S}^{-1/2}\mathbf{V}^\prime$. Then
\[
\mathbf{W}_1 \mathbf{A}_0 \mathbf{W}_2^\prime=
\bigl(\mathbf {W}_1\mathbf{X}_1\bolds{
\Pi}^{1/2}\bigr) \bigl(\mathbf{W}_2\mathbf {X}_2
\bolds{\Pi}^{1/2}\bigr)^\prime= \mathbf{Q}\mathbf{Q}^{-1}
= \eye_r,
\]
where $\eye_r$ denotes the $r\times r$ identity matrix and $\mathbf
{Q}=\mathbf{W}_1\mathbf{X}_1\bolds{\Pi}^{1/2}$.

Moving on, each of $\mathbf{A}_1,\mathbf{A}_2,\ldots,\mathbf
{A}_{\kappa_3}$ has the form
\[
\mathbf{A}_k = \mathbf{X}_1\bolds{\Pi}
\mathbf{D}_k \mathbf {X}_2^\prime, \qquad
\mathbf{D}_k = \operatorname{diag}_k \mathbf{X}_3,
\]
where $\operatorname{diag}_k \mathbf{X}$ denotes the diagonal matrix whose
diagonal equals the $k$th row of matrix $\mathbf{X}$. Applying the
same transformation to $\mathbf{A}_1,\mathbf{A}_2,\ldots,\mathbf
{A}_{\kappa_3}$ yields the collection of $r\times r$ matrices
%
\begin{equation}
\label{eqjointdiag} \mathbf{W}_1 \mathbf{A}_k
\mathbf{W}_2^\prime= \mathbf{Q}\mathbf {D}_k
\mathbf{Q}^{-1}.
\end{equation}
So, the matrices $\lbrace\mathbf{W}_1 \mathbf{A}_k\mathbf
{W}_2^\prime\rbrace$ are diagonalizable in the same basis, namely,
the columns of matrix $\mathbf{Q}$. The associated eigenvalues
$\lbrace\mathbf{D}_k\rbrace$ equal the columns of the matrix
$\mathbf{X}_3$. These eigenvalues are unique up to a joint permutation
of the eigenvectors and eigenvalues provided there exist no $k_1\neq
k_2$ so that the vectors of eigenvalues of $\mathbf{W}_1 \mathbf
{A}_{k_1}\mathbf{W}_2^\prime$ and $\mathbf{W}_1 \mathbf
{A}_{k_2}\mathbf{W}_2^\prime$ are equal [see, e.g., \citet
{DeLathauwerDeMoorVandewalle2004}, Theorem 6.1]. Now, this is
equivalent to demanding that the columns of $\mathbf{X}_3$ are all
distinct. As this is true by assumption, the proof is complete.
\end{pf}

The proof of Theorem \ref{thmid1} shows that access to
lower-dimensional submodels allows to disentangle the scale of the
columns of the $\mathbf{X}_i$ and the weights on the diagonal of
$\bolds{\Pi}$. This is so because the matrix $\bolds{\Pi}$
equally shows up in the lower-dimensional submodels, and so
transforming $\mathbf{A}_k$ to $\mathbf{W}_1\mathbf{A}_k\mathbf
{W}_2^\prime$ absorbs the weights into the joint diagonalizer $\mathbf
{Q}$ in (\ref{eqjointdiag}).

Also note that the dimension of the matrices in (\ref{eqjointdiag})
is $r\times r$, independent of the size of the original matrices
$\mathbf{X}_i$. On the other hand, larger matrices $\mathbf{X}_i$
could be beneficial for identification, as it becomes easier for them
to satisfy the requirement of full column rank.

The full-rank condition that underlies Theorem \ref{thmid1} has a
simple testable implication. Indeed, by (\ref{eqsub}), it implies
that the matrix $\mathbf{A}_0$ has rank $r$. As this matrix is
observable, so is its rank and, hence, our key identifying assumption
is refutable. In applications, this can be done using any of a number
of available rank tests. We refer to \citet{KasaharaShimotsu2014} and
\citet{BonhommeJochmansRobin2014} for practical details on the
implementation of such procedures.

Theorem \ref{thmid1} can be applied to recover the tri-adic
decomposition of $\mathbb{X}$ up to an arbitrary joint permutation
matrix. We present the result in the form of two theorems.

\begin{theorem}[(Vectors)] \label{thmvectors}
If $\mathbf{X}_1$, $\mathbf{X}_2$, and $\mathbf{X}_3$ have full
column rank and for each pair $(i_1,i_2)\in\lbrace i_1,i_2\in\lbrace
1,2,3 \rbrace: i_1<i_2 \rbrace$ $\mathbb{X}_{\lbrace i_1,i_2 \rbrace
}$ is observable, then $\mathbf{X}_1$, $\mathbf{X}_2$, and $\mathbf
{X}_3$ are all identified up to a common permutation of their columns.
\end{theorem}

\begin{theorem}[(Weights)] \label{thmweights}
If $\mathbf{X}_i$ is identified up to a permutation of its columns and
has full column rank, and if $\mathbb{X}_{\lbrace i \rbrace}$ is
observable, then $\bolds{\pi}$ is identified up to the same permutation.
\end{theorem}

\begin{pf}
The one-dimensional submodel $\mathbb{X}_{\lbrace i \rbrace}$ is the vector
\[
\mathbb{X}_{\lbrace i \rbrace} = \mathbf{X}_i \bolds{\pi}.
\]
Given $\mathbf{X}_i$, the one-dimensional submodel yields linear
restrictions on the weight vector $\bolds{\pi}$. Moreover, if
$\mathbf{X}_i$ is known and has maximal column rank, these equations
can be solved for $\bolds{\pi}$, giving
%
\begin{equation}
\label{eqweightsols} \bolds{\pi} = \bigl(\mathbf{X}_i^\prime
\mathbf{X}_i\bigr)^{-1} \mathbf {X}_i^\prime
\mathbb{X}_{\lbrace i \rbrace},
\end{equation}
which is the least-squares coefficient of a regression of $\mathbb
{X}_{\lbrace i \rbrace}$ on the columns of $\mathbf{X}_i$.
\end{pf}

In the supplement, we apply Theorems \ref{thmid1}--\ref{thmweights}
to the finite-mixture model and the hidden Markov model of Section~\ref{secexamples} to obtain constructive proofs of identification.

\section{Estimation by joint approximate diagonalization} \label
{secestimation}
The proof of Theorem \ref{thmid1} shows that the key restrictions
underlying our results take the form of a set of matrices being
simultaneously diagonalizable in the same basis. The problem of joint
matrix diagonalization has recently received considerable attention in
the field of independent component analysis, and
computationally-efficient algorithms for it have been developed; see
\citet{FuGao2006},
\citeauthor{IferroudjeneAbedMeraimBelouchrani2009} (\citeyear{IferroudjeneAbedMeraimBelouchrani2009,IferroudjeneAbed-MeraimBelouchrani2010})
and \citeauthor{LucianiAlbera2010} (\citeyear{LucianiAlbera2010,LucianiAlbera2014}). Such algorithms can be
exploited here to construct easy-to-implement nonparametric estimators
of multivariate latent-structure models.

Thus, we propose estimating the latent-structure model in (\ref
{eqqad}) as follows. Given an estimate of the array $\mathbb{X}$ and
of its lower-dimensional submodels, first estimate all $\mathbf
{x}_{ij}$ by solving a sample version of the joint diagonalization
problem in~(\ref{eqjointdiag}), possibly after unfolding if $q>3$.
Next, back out the weights $\pi_1,\pi_2,\ldots,\pi_r$ by solving
the sample analog of the minimum-distance problem in~(\ref
{eqweightsols}). Asymptotic theory for this second step follows
readily by the delta method. If desired, a consistent labelling can be
recovered based on the proof of Theorem \ref{thmvectors} (see the
supplementary material).

\subsection{Estimator}
Consider a generic situation in which a set of $\kappa r\times r$
matrices $\mathbf{C}_{1},\mathbf{C}_{2},\ldots,\mathbf{C}_{\kappa
}$ can be jointly diagonalized by an $r\times r$ invertible matrix
$\mathbf{Q}_{0}$, that is,
%
\begin{equation}
\label{eqsystem1} \mathbf{C}_k = \mathbf{Q}_0
\mathbf{D}_k\mathbf{Q}_0^{-1},
\end{equation}
for diagonal matrices $\mathbf{D}_{1},\mathbf{D}_{2},\ldots,\mathbf
{D}_{\kappa}$. Knowledge of the joint eigenvectors implies knowledge
of the eigenvalues as
%
\begin{equation}
\label{eqsystem2} \mathbf{D}_k = \mathbf{Q}_0^{-1}
\mathbf{C}_k \mathbf{Q}_0.
\end{equation}
The matrix $\mathbf{Q}_0$ is not unique. Moreover, let $\off\mathbf
{Q} = \mathbf{Q}-\operatorname{diag}\mathbf{Q}$ and let $\llVert
\mathbf{Q}\rrVert _F = \sqrt{\tr(\mathbf{Q}^\prime\mathbf
{Q})}$ denote the Frobenius norm. Then any solution to the
least-squares problem
%
\begin{equation}
\label{eqols} \min_{\mathbf{Q}} \sum_{k=1}^\kappa
\bigl\llVert \off\bigl(\mathbf {Q}^{-1}\mathbf{C}_k
\mathbf{Q}\bigr) \bigr\rrVert _F^2
\end{equation}
is a joint diagonalizer in the sense of (\ref{eqsystem1}). Each of
these delivers the same set of eigenvalues in (\ref{eqsystem2}) (up
to a joint permutation).

The statistical problem of interest in this section is to perform
inference on the $\mathbf{D}_1,\mathbf{D}_2,\ldots,\mathbf
{D}_\kappa$ when\vspace*{1pt} we only observe noisy versions of the input matrices
$\mathbf{C}_1,\mathbf{C}_2,\ldots,\mathbf{C}_\kappa$, say
$\widehat{\mathbf{C}}_1,\widehat{\mathbf{C}}_2,\ldots,\widehat
{\mathbf{C}}_\kappa$. The sampling noise in the $\widehat{\mathbf
{C}}_{k}$ prevents them from sharing the same set of eigenvectors.
Indeed, in general, there does not exist a $\mathbf{Q}$ such that
$\mathbf{Q}^{-1}\widehat{\mathbf{C}}_{k} \mathbf{Q}$
will be exactly diagonal for all $k$. For this, the least-squares
formulation in (\ref{eqsystem2})--(\ref{eqols}) is important as it
readily suggests using, say $\widehat{\mathbf{Q}}$, any solution to
%
\begin{equation}
\label{eqolssample} \min_{\mathbf{Q}\in\mathscr{Q}} \sum_{k=1}^\kappa
\bigl\llVert \off \bigl(\mathbf{Q}^{-1}\widehat{\mathbf{C}}_k
\mathbf{Q}\bigr) \bigr\rrVert _F^2,
\end{equation}
where\vspace*{2pt} $\mathscr{Q}$ is an appropriately-specified space of matrices to
search over; see below. The estimator $\widehat{\mathbf{Q}}$ is that
matrix that makes all these matrices as diagonal as possible, in the
sense of minimizing the sum of their squared off-diagonal entries. It
is thus appropriate\vspace*{2pt} to call the estimator $\widehat{\mathbf{Q}}$ the
joint approximate-diagonalizer of $\widehat{\mathbf{C}}_1,\widehat
{\mathbf{C}}_2,\ldots,\widehat{\mathbf{C}}_\kappa$. An estimator
of the $\mathbf{D}_k$ (up to a joint permutation of their eigenvalues)
then is
%
\begin{equation}
\label{eqDhat} \widehat{\mathbf{D}}_k = \operatorname{diag}\bigl(
\widehat{\mathbf {Q}}^{-1}\widehat{\mathbf{C}}_k \widehat{
\mathbf{Q}}\bigr).
\end{equation}
Distribution theory for this estimator is not available, however, and
so we provide it here. Throughout, we work under the convention that
estimates are computed from a sample of size $n$.

\subsection{Asymptotic theory}
For our problem to be well defined, we assume that the matrix of joint
eigenvectors is bounded. In (\ref{eqolssample}), we may therefore
restrict attention to the set of $r\times r$ matrices $\mathbf
{Q}=(\mathbf{q}_1,\mathbf{q}_2,\ldots,\mathbf{q}_r)$ defined as
\[
\mathscr{Q} = \bigl\lbrace\mathbf{Q}: \det\mathbf{Q}=1, \llVert
\mathbf{q}_j\rrVert _F= c\mbox{ for } j=1,2,
\ldots,r\mbox{ and }c\leq m \bigr\rbrace
\]
for some $m\in(0,\infty)$. The restrictions on the determinant and
the column norms are without loss of generality and only reduce the
space of matrices to be searched over when solving (\ref
{eqolssample}). Let $\mathbf{Q}_*$ be any solution to (\ref
{eqols}) on $\mathscr{Q}$ and let $\mathscr{Q}_0\subset\mathscr
{Q}$ be the set of all matrices $\mathbf{Q}_*\bolds{\Delta} \bolds
{\Theta}$ for permutation matrices $\bolds{\Delta}$ and diagonal
matrices $\bolds{\Theta}$ whose diagonal entries are equal to $1$
and $-1$ and have $\det\bolds{\Theta}=1$. Then $\mathscr{Q}_0$ is
the set of solutions to (\ref{eqols}) on $\mathscr{Q}$.

Construct the $r\times r\kappa$ matrix $\mathbf{C}=(\mathbf
{C}_1,\mathbf{C}_2,\ldots,\mathbf{C}_\kappa)$ by concatenation and
define $\widehat{\mathbf{C}}$ similarly.

\begin{theorem}[(Consistency)] \label{thmconsistency}
If\vspace*{2pt} the set $\mathscr{Q}_0$ belongs to the interior of $\mathscr{Q}$,
$\widehat{\mathbf{C}}=\mathbf{C}+o_p(1)$, and $\widehat{\mathbf
{Q}}\in\mathscr{Q}$ satisfies
\[
\sum_{k=1}^\kappa\bigl\llVert \off\bigl(
\widehat{\mathbf {Q}}^{-1}\widehat{\mathbf{C}}_k \widehat{
\mathbf{Q}}\bigr) \bigr\rrVert _F^2 = \min
_{\mathbf{Q}\in\mathscr{Q}} \Biggl\lbrace \sum_{k=1}^\kappa
\bigl\llVert \off\bigl(\mathbf{Q}^{-1}\widehat{\mathbf
{C}}_k \mathbf{Q}\bigr) \bigr\rrVert _F^2
\Biggr\rbrace+ o_p(1),
\]
then $\lim_{n\rightarrow\infty}\Pr(\widehat{\mathbf{Q}}\in
\mathscr{O}) =1$ for any open subset $\mathscr{O}$ of $\mathscr{Q}$
containing $\mathscr{Q}_0$.
\end{theorem}

Each $\mathbf{Q}\in\mathscr{Q}_0$ has associated with it a
permutation matrix $\bolds{\Delta}$ and a diagonal matrix $\bolds
{\Theta}$ as just defined so that $\mathbf{Q}=\mathbf{Q}_* \bolds
{\Delta} \bolds{\Theta}$. Theorem \ref{thmconsistency} states
that (up to a subsequence) we have that
$
\widehat{\mathbf{Q}} \overset{p}{\rightarrow} \mathbf{Q}_* \bolds
{\Delta}_0 \bolds{\Theta}_0
$
for well-defined $\bolds{\Delta}_0$ and $\bolds{\Theta}_0$. We
may then set $\mathbf{Q}_0 = \mathbf{Q}_*\bolds{\Delta}_0 \bolds
{\Theta}_0$ in (\ref{eqsystem1}). It then equally follows that
\[
\widehat{\mathbf{D}}_k \overset{p} {\rightarrow}
\mathbf{D}_k = \bolds{\Delta}_0^\prime
\mathbf{D}_k^* \bolds{\Delta}_0,
\]
where $\mathbf{D}_k$ is as in (\ref{eqsystem2}) and $\mathbf{D}_k^*
= \mathbf{Q}_*^{-1}\mathbf{C}_k \mathbf{Q}_*$, both of which are
equal up to a permutation. Thus, the consistency of the eigenvalues (up
to a joint permutation) follows from the consistency of the estimator
of the input matrices $\mathbf{C}$.

To provide distribution theory, let
\[
\mathbf{D}_{k_1}\ominus\mathbf{D}_{k_2} = (
\mathbf{D}_{k_1} \otimes \eye_{\dim \mathbf{D}_{k_2}}) - (\eye_{\dim\mathbf
{D}_{k_1}}
\otimes\mathbf{D}_{k_2})
\]
denote the Kronecker difference between the square matrices $\mathbf
{D}_{k_1}$ and $\mathbf{D}_{k_{2}}$. Construct the $r^2\times
r^2\kappa$ matrix
\[
\mathbf{T} = \bigl((\mathbf{D}_{1}\ominus\mathbf{D}_{1}),(
\mathbf{D}_{2}\ominus \mathbf{D}_{2}),\ldots,(
\mathbf{D}_{\kappa}\ominus\mathbf {D}_{\kappa}) \bigr)
\]
by concatenation and let
\[
\mathbf{G} = (\eye_r \otimes \mathbf{Q}_0) \Biggl(\sum
_{k=1}^\kappa (\mathbf{D}_{k}\ominus
\mathbf{D}_{k})^2 \Biggr)^+ \mathbf{T} \bigl(
\eye_\kappa\otimes \mathbf{Q}_0^\prime \otimes
\mathbf{Q}_0^{-1}\bigr),
\]
where $\mathbf{Q}^+$ is the Moore--Penrose pseudo inverse of $\mathbf
{Q}$. Theorem \ref{thmnormalityQ} contains distribution theory for
our estimator of the matrix of joint eigenvectors $\widehat{\mathbf
{Q}}$ in (\ref{eqolssample}).

\begin{theorem}[(Asymptotic distribution)]\label{thmnormalityQ}
If $\llVert \widehat{\mathbf{C}}-\mathbf{C}\rrVert _F =
O_p(n^{-1/2})$, then
\[
\sqrt{n}\vc(\widehat{\mathbf{Q}}-\mathbf{Q}_0) = \mathbf{G}
\sqrt{n}\vc(\widehat{\mathbf{C}}-\mathbf{C}) + o_p(1)
\]
as $n\rightarrow\infty$.
\end{theorem}

If, further, $\sqrt{n}  \operatorname{vec} (\widehat{\mathbf{C}}-\mathbf
{C})\overset{d}{\rightarrow}\mathcal{N}(\boldsymbol{0},\mathbf
{V})$ for some covariance matrix $\mathbf{V}$, Theorem~\ref
{thmnormalityQ} implies that
\[
\sqrt{n} \operatorname{vec} (\widehat{\mathbf{Q}} - \mathbf {Q}_0)
\overset{d} {\rightarrow}\mathcal{N}\bigl(\boldsymbol{0},\mathbf {G} \mathbf{V}
\mathbf{G}^\prime\bigr)
\]
as $n\rightarrow\infty$. In our context, $\sqrt{n}$-consistency and
asymptotic normality of the input matrices is not a strong requirement.
Indeed, the proof of Theorem \ref{thmid1} showed that the input
matrices are of the form $\mathbf{C}_k=\mathbf{W}_1 \mathbf{A}_k
\mathbf{W}_2^\prime$, where $\mathbf{W}_1$ and $\mathbf{W}_2$
follow from a singular-value decomposition of $\mathbf{A}_0$. An
estimator of $\mathbf{C}_k$ can thus be constructed using a sample
analog of $\mathbf{A}_0$ to estimate $\mathbf{W}_1$ and $\mathbf
{W}_2$, together with a sample analog of $\mathbf{A}_k$. If the
estimators of $\mathbf{A}_0$ and $\mathbf{A}_k$ are $\sqrt
{n}$-consistent and asymptotically normal and all nonzero singular
values of $\mathbf{A}_0$ are simple, then $\sqrt{n}  \operatorname{vec}
(\widehat{\mathbf{C}}-\mathbf{C})\overset{d}{\rightarrow}\mathcal
{N}(\boldsymbol{0},\mathbf{V})$ holds. A detailed derivation of
$\mathbf{V}$ is readily obtained from the argument on the estimation
of eigen-decompositions of normal matrices in the supplementary
material to
\citeauthor{BonhommeJochmansRobin2014} [(\citeyear{BonhommeJochmansRobin2014}), Lemma S.2].\vspace*{1pt}

We next present the asymptotic behavior of
$\widehat{\mathbf{D}}=(\widehat{\mathbf{D}}_1,\widehat{\mathbf
{D}}_2,\ldots,\widehat{\mathbf{D}}_\kappa)$, our estimator of the
eigenvalues $\mathbf{D}=(\mathbf{D}_1,\mathbf{D}_2,\ldots,\mathbf
{D}_\kappa)$. To state it, let $\mathbf{S}_r = \operatorname{diag}(\vc
\eye_r)$ be an $r^2\times r^2$ selection matrix; note that $\mathbf
{S}_r\vc\mathbf{Q} = \vc(\operatorname{diag}\mathbf{Q})$. Let
\[
\mathbf{H} = (\eye_\kappa\otimes \mathbf{S}_r) \bigl(
\eye_\kappa\otimes \mathbf{Q}_0^\prime \otimes\mathbf
{Q}_0^{-1}\bigr).
\]
Theorem \ref{thmnormality} follows.

\begin{theorem}[(Asymptotic distribution)]\label{thmnormality}
If $\llVert \widehat{\mathbf{C}}-\mathbf{C}\rrVert _F =
O_p(n^{-1/2})$, then
\[
\sqrt{n} \operatorname{vec} (\widehat{\mathbf{D}} - \mathbf{D}) = \mathbf{H} \sqrt{n}
\operatorname{vec} (\widehat{\mathbf{C}}- \mathbf{C}) +o_p(1)
\]
as $n\rightarrow\infty$.
\end{theorem}

Again, if $\sqrt{n}  \operatorname{vec} (\widehat{\mathbf{C}}-\mathbf
{C})\overset{d}{\rightarrow}\mathcal{N}(\boldsymbol{0},\mathbf
{V})$, then
\[
\sqrt{n} \operatorname{vec} (\widehat{\mathbf{D}}-\mathbf{D})\overset {d} {
\rightarrow}\mathcal{N}\bigl(\boldsymbol{0},\mathbf{H} \mathbf{V}
\mathbf{H}^\prime\bigr)
\]
as $n\rightarrow\infty$.

\section{Application to density estimation} \label{secmixtures}
With discrete outcomes, both the finite-mixture model in (\ref
{eqdiscretemixture}) and the hidden Markov model in (\ref{eqhmm})
are finite dimensional. Further, the matrices to be simultaneously
diagonalized are contingency tables. These tables can be estimated by
simple empirical cell probabilities and are $\sqrt{n}$-consistent and
asymptotically normal. Hence, the theory on the asymptotic behavior of
the eigenvalues from the previous section (i.e., Theorem \ref
{thmnormality}) can directly be applied to deduce the large-sample
behavior of the parameter estimates.

With continuous outcomes, as in (\ref{eqcontinuousmixture}), the main
parameters of the model are density functions. Such an
infinite-dimensional problem is not directly covered by the arguments
from the previous section. Nonetheless, we will show that Theorem~\ref
{thmnormalityQ} can be used to obtain density estimators with standard
asymptotic properties.

\subsection{Estimator}
We provide convergence rates and distribution theory for series
estimators based on (\ref{eqfouriersystem}). By the results of
Section~\ref{subsechmm}, this also covers the estimation of
emission densities in a hidden Markov model with continuous outcome
variables. Recall from above that the projections
\[
\proj_{\kappa_i} f_{ij} = \bolds{\varphi}_{\kappa_i}^\prime
\mathbf{b}_{ij}
\]
yield the multilinear restrictions
\[
\mathbb{B}=E\Biggl[\bigotimes_{i=1}^q
\bolds{\varphi}_{\kappa
_i}(Y_i)\rho(Y_i)\Biggr] =
\sum_{j=1}^r
\pi_j \bigotimes_{i=1}^q E
\bigl[\bolds{\varphi}_{\kappa
_i}(Y_i)\rho(Y_i)\mid
Z=j\bigr] = \sum_{j=1}^r
\pi_j \bigotimes_{i=1}^q
\mathbf{b}_{ij},
\]
where $\bolds{\varphi}_{\kappa_i}$ is the vector containing the
$\kappa_i$ leading polynomials from the orthogonal system $\lbrace
\varphi_k, k>0\rbrace$. As we will show, for fixed $\kappa_1,\kappa
_2,\ldots,\kappa_q$, the array $\mathbb{B}$ provides sufficient
information for nonparametric identification of Fourier coefficients
through the associated joint diagonalizer. Moreover, in the asymptotic
analysis, $\kappa_1,\kappa_2,\ldots,\kappa_q$ are all held fixed.

For the purpose of this section, we may fix attention to a given index
$i$. By unfolding $\mathbb{B}$ toward direction $i$, we obtain the
(equivalent) three-way array
\[
\mathbb{B}_i = E \bigl[\bolds{\phi}^{\mathcal{Q}_1}\otimes\bolds {
\phi}^{\mathcal{Q}_2} \otimes\bolds{\varphi}_{\kappa
_i}(Y_{i})
\rho(Y_{i}) \bigr],
\]
where $\mathcal{Q}_1$ and $\mathcal{Q}_2$ partition the index set
$\lbrace1,2,\ldots,q \rbrace\setminus\lbrace i \rbrace$ (see
Section~\ref{secid}) and we have introduced the notational shorthand
\[
\bolds{\phi}^{\mathcal{Q}} = \bigodot_{i^\prime\in\mathcal{Q}} \bolds{
\varphi}_{\kappa
_{i^\prime}}(Y_{{i^\prime}})\rho(Y_{{i^\prime}}).
\]
The array $\mathbb{B}_i$ can be analyzed using our diagonalization
approach. Following the notation from the proof of Theorem~\ref
{thmid1}, the two-dimensional submodel associated with $\mathbb{B}_i$
is the matrix
\[
\mathbf{A}_0 = E \bigl[ \bolds{\phi}^{\mathcal{Q}_1}\otimes\bolds
{\phi}^{\mathcal{Q}_2} \bigr],
\]
while the array $\mathbb{B}_i$ itself consists of the first $\kappa
_i$ matrices of the set $\lbrace\mathbf{A}_k, k>0\rbrace$, where
\[
\mathbf{A}_k = E \bigl[ \bigl( \bolds{\phi}^{\mathcal{Q}_1}\otimes
\bolds{\phi }^{\mathcal{Q}_2} \bigr) \varphi_k(Y_{i})
\rho(Y_{i}) \bigr].
\]
All these matrices are of dimension $\prod_{i_1\in\mathcal
{Q}_1}\kappa_{i_1}\times\prod_{i_2\in\mathcal{Q}_2}\kappa_{i_2}$.
A singular-value decomposition of $\mathbf{A}_0$ provides matrices
$\mathbf{W}_1$ and $\mathbf{W}_2$ so that the $\kappa_i$ matrices
$
\mathbf{W}_1\mathbf{A}_k\mathbf{W}_2^\prime
$
are jointly diagonalizable by, say, $\mathbf{Q}$. From the proof of
Theorem~\ref{thmid1}, the matrix $\mathbf{Q}$ is unique (up to the
usual normalizations on the sign and norm of its columns and a joint
permutation of the columns, as discussed before) as soon as the
conditions in Theorem~\ref{thmid1} are satisfied.

Given $\mathbf{Q}$, we can compute
\[
\mathbf{Q}^{-1}\bigl(\mathbf{W}_1\mathbf{A}_k
\mathbf{W}_2^\prime\bigr) \mathbf{Q} = \operatorname{diag}
(b_{i1k},b_{i2k},\ldots,b_{irk}),
\]
where, recall,
$
b_{ijk} = E[\varphi_k(Y_i)\rho(Y_i) \mid Z=j]
$
for any integer $k$ (including those $k$ that exceed $\kappa_i$).
Equivalently, the $k$th Fourier coefficient of $f_{ij}$ can be written as
%
\begin{equation}
\label{eqcoeff} b_{ijk} = \mathbf{e}_j^\prime
\bigl( \mathbf{Q}^{-1}\bigl(\mathbf{W}_1
\mathbf{A}_k\mathbf{W}_2^\prime\bigr) \mathbf{Q}
\bigr) \mathbf{e}_j,
\end{equation}
where $\mathbf{e}_j$ is the $r\times1$ selection vector whose $j$th
entry is equal to one and its other entries are all equal to zero.

Our orthogonal-series estimator of $f_{ij}$ is based on sample analogs
of the $b_{ijk}$ in~(\ref{eqcoeff}). We estimate the array $\mathbb
{B}$ as
\[
\widehat{\mathbb{B}} = n^{-1} \sum_{m=1}^n
\bigotimes_{i=1}^q \bolds{
\varphi}_{\kappa_i}(Y_{im}) \rho(Y_{im}),
\]
where $\lbrace Y_{1m},Y_{2m},\ldots,Y_{qm} \rbrace_{m=1}^n$ is a
size-$n$ sample drawn at random from the mixture model. From this we
estimate $b_{ijk}$ for any $k$ as
\[
\hat{b}_{ijk} = \mathbf{e}_j^\prime \bigl(
\widehat{\mathbf{Q}}^{-1}\bigl(\widehat{\mathbf{W}}_1
\widehat{\mathbf {A}}_k\widehat{\mathbf{W}}_2^\prime
\bigr) \widehat{\mathbf{Q}} \bigr) \mathbf{e}_j = n^{-1}
\sum_{m=1}^n \mathbf{e}_j^\prime
\widehat{\bolds{\Omega}}_m\mathbf{e}_j
\varphi_k(Y_{im}) \rho(Y_{im}),
\]
using obvious notation to denote sample counterparts in the first
expression and introducing the matrix
\[
\widehat{\bolds{\Omega}}_m = \widehat{\mathbf{Q}}^{-1}
\bigl(\widehat{\mathbf{W}}_1 \bigl(\bolds{\phi}_m^{\mathcal{Q}_1}
\otimes\bolds{\phi }_{m}^{\mathcal{Q}_2} \bigr) \widehat{
\mathbf{W}}_2^\prime\bigr) \widehat{\mathbf{Q}}
\]
in the second expression; here, we let
$
\bolds{\phi}_m^{\mathcal{Q}}
=
\bigodot_{i^\prime\in\mathcal{Q}} \bolds{\varphi}_{\kappa
_{i^\prime}}(Y_{i^\prime m})\rho(Y_{i^\prime m})$.
The associated orthogonal-series estimator of $f_{ij}(y)$ for some
chosen integer $\varkappa$ is
%
\begin{eqnarray}
\label{eqdensityhat} \hat{f}_{ij}(y) &=& \sum_{k=1}^{\varkappa}
\hat{b}_{ijk} \varphi_{k}(y)
\nonumber\\[-8pt]\\[-8pt]\nonumber
& =& n^{-1} \sum_{m=1}^n
\mathbf{e}_j^\prime\widehat{\bolds{\Omega}}_m
\mathbf{e}_j \sum_{k=1}^{\varkappa}
\varphi_k(Y_{im})\varphi_k(y)
\rho(Y_{im}).
\end{eqnarray}
Note that, in the absence of $\mathbf{e}_j^\prime\widehat{\bolds
{\Omega}}_m\mathbf{e}_j$, this expression collapses to a standard
series estimator of the marginal density of $Y_i$. Hence,\vspace*{2pt} the term
$\mathbf{e}_j^\prime\widehat{\bolds{\Omega}}_m\mathbf{e}_j$ can
be understood as a weight that transforms this estimator into one of
the conditional density of $Y_i$ given $Z=j$. Equation (\ref
{eqdensityhat}) generalizes the kernel estimator of \citet{BonhommeJochmansRobin2014}. The term $\mathbf{e}_j^\prime\widehat
{\bolds{\Omega}}_m\mathbf{e}_j$ plays\vspace*{1pt} the same role as the
posterior classification probability (normalized to sum up to one
across observations) in the EM algorithm as well as in its
nonparametric version [\citet{LevineHunterChauveau2011}, equations
(15)--(17)]. A computational advantage here
is that the series estimator is available in closed form once $\mathbf
{e}_j^\prime\widehat{\bolds{\Omega}}_m\mathbf{e}_j$ has been
computed while EM requires iterative computation of density estimates
and classification probabilities until convergence.

A natural way of choosing the number of series terms in (\ref
{eqdensityhat}) would be by minimizing the squared $L_\rho^2$-loss,
\[
\llVert \hat{f}_{ij}-f_{ij} \rrVert _2^2,
\]
as a function of $\varkappa$. In the supplement we show that an
empirical counterpart of this criterion (up to terms that do not
involve $\varkappa$) is
\[
\sum_{k=1}^\varkappa\hat{b}_{ijk}^2
- \frac{2n^{-1}}{n-1} \sum_{m=1}^n \sum
_{o\neq m} \mathbf{e}_j^\prime
\widehat{\bolds{\Omega}}_m\mathbf{e}_j
\mathbf{e}_j^\prime\widehat{\bolds{\Omega}}_o
\mathbf{e}_j \sum_{k=1}^{\varkappa}
\varphi_k(Y_{io})\varphi_k(Y_{im})
\rho(Y_{io}) \rho(Y_{im}).
\]
Apart from the weight functions, this is the usual cross-validation
objective for orthogonal-series estimators [\citet{Hall1987}].

Before turning to the statistical properties of $\hat{f}_{ij}$ we note
that, although we maintain a hard thresholding procedure in (\ref
{eqdensityhat}), our approach can equally be combined with other
popular smoothing policies that shrink the impact of higher-order
Fourier coefficients; see
\citeauthor{Efromovich1999} [(\citeyear{Efromovich1999}), Chapter~3] for a discussion on such policies.

\subsection{Asymptotic theory}
Under mild conditions, the series estimator in (\ref{eqdensityhat})
exhibits standard large-sample behavior. The precise conditions depend
on the choice of orthogonal system, that is, $\lbrace\varphi_k,k>0
\rbrace$. We give two sets of conditions that cover the most popular choices.

When the component densities are supported on compact intervals, we can
restrict attention to $[-1,1]$ without loss of generality; translation
to generic compact sets is straightforward. In this case, we will allow
for polynomial systems that satisfy the following general requirements.
Here and later, we let $\llVert \cdot\rrVert _\infty$ denote
the supremum norm.
\begin{longlist}[A.1]
\item[A.1] The sequence $\lbrace\varphi_k,k> 0 \rbrace$ is dominated by a
function $\psi$, which is continuous on $(-1,1)$ and positive almost
everywhere on $[-1,1]$. $\rho$, $\psi\rho$, and $\psi^2\rho$ are
integrable, and there exists a sequence of constants $\lbrace\zeta
_\varkappa, \varkappa> 0\rbrace$ so that $\llVert \sqrt{\bolds
{\varphi}_\varkappa^\prime\bolds{\varphi}_\varkappa} \rrVert _\infty\leq\zeta_\varkappa$.
\end{longlist}
These conditions are rather weak. They are satisfied for the popular
class of Jacobi polynomials, for example, which includes Chebyshev
polynomials of the first kind, Chebyshev polynomials of the second
kind, and Legendre polynomials.

In this case, we will need the following regularity from the component
densities.
\begin{longlist}[A.2]
\item[A.2] The $(\psi\rho)^{4}f_{ij}$ are integrable.
\end{longlist}
The weaker requirement that the $(\psi\rho)^{2}f_{ij}$ are integrable
will suffice to obtain the convergence rates in Theorem~\ref
{thmrates} below, but \textup{A.2} will be needed to obtain the
pointwise asymptotic-normality result in Theorem~\ref{thmpointwise}.

When the component densities are supported on the whole real line, we
will take $\lbrace\varphi_k,k> 0 \rbrace$ to be the orthonormalized
system of Hermite functions.
\begin{longlist}[B.1]
\item[B.1] The sequence $\lbrace\varphi_k,k> 0 \rbrace$ has members
\[
\varphi_k(y) = 2^{-(k-1)/2} \bigl((k-1){!}\bigr)^{-1/2}
\pi^{-1/4} e^{-y^2/2} h_{k-1}(y),
\]
where $\lbrace h_k, k\geq0 \rbrace$ is the system of the Hermite
polynomials, in which case $\llVert \sqrt{\bolds{\varphi
}_\varkappa^\prime\bolds{\varphi}_\varkappa} \rrVert _\infty
\leq\zeta_\varkappa$ for $\zeta_\varkappa\propto\sqrt{\varkappa}$.
\end{longlist}

We will also impose the following regularity and smoothness conditions.
\begin{longlist}[C.3]
\item[C.1] The $f_{ij}$ are continuous.

\item[C.2] $\llVert \proj_{\varkappa}f_{ij}-f_{ij}\rrVert
_{\infty}=O(\varkappa^{-\beta})$ for some constant $\beta\geq1$.

\item[C.3] The singular values of $\mathbf{A}_0$ are all simple.
\end{longlist}
Convergence in $L_\rho^2$-norm implies that $\lim_{\varkappa
\rightarrow\infty} \sum_{k=1}^\varkappa b_{ijk}^2$ is finite, and so
that the Fourier coefficient associated with $\varphi_k$ shrinks to
zero as $k\rightarrow\infty$. The constant $\beta$ is a measure of
how fast the Fourier coefficients shrink. In general, $\beta$ is
larger the smoother the underlying function that is being approximated.
Simplicity of the singular values of $\mathbf{A}_0$ holds generically
and is used here to ensure that the matrices $\mathbf{W}_1,\mathbf
{W}_2$ are continuous transformations of $\mathbf{A}_0$. This is a
technical requirement used to derive the convergence rates of their
plug-in estimators.

Under these assumptions, we obtain standard integrated squared-error
and uniform convergence rates.

\begin{theorem}[(Convergence rates)] \label{thmrates}
Let either \textup{A.1}--\textup{A.2} and \textup{C.1}--\textup{C.3} or \textup{B.1} and \textup{C.1}--\textup{C.3} hold. Then
\[
\llVert \hat{f}_{ij} -f_{ij} \rrVert _2^2
= O_p\bigl(\varkappa/n+\varkappa^{-2\beta}\bigr), \qquad
\llVert \hat{f}_{ij} -f_{ij} \rrVert _\infty=
O_p\bigl(\zeta_{\varkappa}\sqrt{\varkappa/n}+
\varkappa^{-\beta}\bigr),
\]
for all $i,j$.
\end{theorem}

The rates in Theorem~\ref{thmrates} equal the conventional univariate
rates of series estimators; see, for example, \citet{Newey1997}. Thus,
the fact that $Z$ is latent does not affect the convergence speed of
the density estimates.

To present distribution theory for the orthogonal-series estimator at a
fixed point~$y$, let
\[
\hat{\sigma}_{ij}(y) = \sqrt{ n^{-1} \sum
_{m=1}^n \Biggl( \mathbf{e}_j^\prime
\widehat{\bolds{\Omega}}_m\mathbf{e}_j \sum
_{k=1}^\varkappa\varphi_k(Y_{im})
\varphi_k(y) \rho(Y_{im}) - \hat{f}_{ij}(y)
\Biggr)^2 },
\]
which is a sample standard deviation, and denote $f_i = \sum_{j=1}^r
\pi_j   f_{ij}$ in the following theorem.

\begin{theorem}[(Asymptotic distribution)] \label{thmpointwise}
Suppose that $n,\varkappa\rightarrow\infty$ so that $\varkappa
^2/n\rightarrow0$ and $n\varkappa^{-2\beta}\rightarrow0$. Then
\[
\frac{\hat{f}_{ij}(y)-f_{ij}(y)}{\hat{\sigma}_{ij}(y)/\sqrt
{n}}\overset{d} {\rightarrow}\mathcal{N}(0,1),
\]
for each $y\in\mathscr{Y}$ that lies in an interval on which $f_i$ is
of bounded variation.
\end{theorem}

Under \textup{A.1}--\textup{A.2}, $\hat{\sigma}_{ij}(y)$ grows like
$\llVert \bolds{\varphi}_{\varkappa}(y) \rrVert _F$, and
this depends on the polynomial system used. Because \textup{A.1}
states that $\llVert \sqrt{\boldsymbol{\varphi}_{\varkappa
}^\prime\boldsymbol{\varphi}_{\varkappa}} \rrVert _\infty
=O(\zeta_{\varkappa})$, a weak bound on the convergence rate that
holds for all $y$ is $O_p(\zeta_{\varkappa}/\sqrt{n})$. With
Legendre polynomials, for example, the orthogonal-series estimator has
a variance of order $\varkappa/n$, which is the same as that of an
estimator based on a random\vspace*{1pt} sample from $f_{ij}$ [\citet
{Hall1987}]. Likewise, under \textup{B.1} we have that $\hat{\sigma
}_{ij}(y)$ grows like ${\varkappa^{1/4}}$ and so the variance of the
estimator is of the order ${\sqrt{\varkappa}/n}$. This is again the
standard convergence rate for conventional Hermite series estimators
[\citet{Liebsher1990}].

\section{Monte Carlo illustrations}
We evaluated the performance of the ortho\-gonal-series estimator via
simulation. We report root mean integrated squared error (RMISE)
calculations for designs taken from \citet{LevineHunterChauveau2011}.
This allows us to compare our estimator to the EM-like approaches
proposed in the literature. We also investigate the accuracy of the
pointwise asymptotic approximation of the density estimator in Theorem
\ref{thmpointwise} in a Monte Carlo experiment based on a hidden
Markov model. Throughout this section, we use Hermite polynomials as
basis functions, set $\kappa_i=10$ for all $i$, and use the
cross-validation technique introduced above to select the number of
series terms. Joint approximate diagonalization was done using the
algorithm of \citeauthor{LucianiAlbera2010} (\citeyear{LucianiAlbera2010,LucianiAlbera2014}). We also
computed the estimator using the algorithms of \citet{FuGao2006} and
\citeauthor{IferroudjeneAbedMeraimBelouchrani2009} (\citeyear{IferroudjeneAbedMeraimBelouchrani2009,IferroudjeneAbed-MeraimBelouchrani2010})
and found very similar results to the ones reported below.

\subsection{RMISE comparisons} \label{subsecRMISE}
We evaluate the RMISE of the estimator $\hat{f}_{ij}$,
\[
\sqrt{E\llVert \hat{f}_{ij}-f_{ij}\rrVert
_2^2},
\]
as approximated by $500$ Monte Carlo replications. The first set of
designs involves mixtures of normals, where
\[
f_{ij}(y) = \phi (y-\mu_{ij} ).
\]
The second set of designs deals with mixtures of central and
noncentral $t$-distri\-butions, that is,
\[
f_{ij}(y) = t_{10}(y;\mu_{ij}),
\]
where we let $t_d(y;\mu)$ denote a $t$-distribution with $d$ degrees
of freedom and noncentrality parameter $\mu$. We set $q=3$, $r=2$, so
the data is drawn from a three-variate two-component mixture. The
parameters of the component densities are
set to $(\mu_{11},\mu_{21},\mu_{31})=(0,0,0)$ for the first
component and $(\mu_{12},\mu_{22},\mu_{32})=(3,4,5)$ for the second
component. We consider various choices for the mixing proportions
$\bolds{\pi}=(\pi_1,\pi_2)^\prime$.

\begin{figure}

\includegraphics{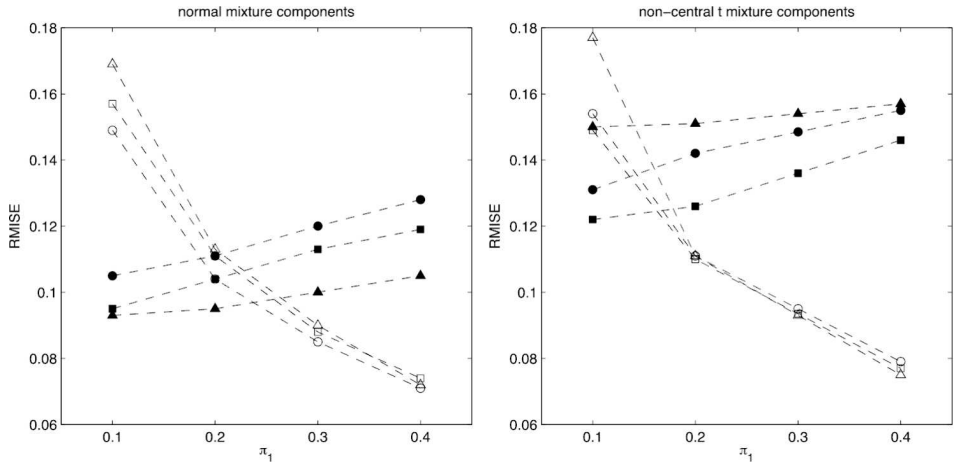}

\caption{RMISE of the orthogonal-series density estimator.}
\label{figmcrmise}
\end{figure}

Figure~\ref{figmcrmise} plots the RMISE as a function of the mixing
proportion $\pi_1$ for samples of size $n=500$. The results for the
first and second component for each outcome variable are labelled
consecutively as $\circ,\square,\triangle$ and as $\bullet$, $\blacksquare$, $\blacktriangle$, respectively.

The patterns of the RMISE are comparable to those for the EM-like
estimators in
\citeauthor{LevineHunterChauveau2011} [(\citeyear{LevineHunterChauveau2011}), Figure~1], although the
magnitudes are larger here. The latter observation agrees with the
intuition that joint estimation of classification probabilities and
component densities (as in EM) should be more efficient than sequential
estimation (as here). However, a precise comparison between the methods
is complicated by the fact that the EM approaches are kernel based
while we work with orthogonal series, and because the tuning parameters
(the bandwidths for EM and the number of series terms here) were
selected in a different manner.

Our least-squares estimator of the mixing proportions was also
evaluated in these designs and was found to perform well. The Monte
Carlo results are provided in the supplementary material [\citet{BonhommeJochmansRobin2014b}].

\subsection{Inference in a hidden Markov model} \label{subsecHMM}
We next consider inference in a hidden Markov model with $r=2$ latent
states and $q=3$ outcome variables. The latent Markov chain has
transition matrix and stationary distribution equal to
\[
\mathbf{K} = \pmatrix{ 0.8 & 0.2
\cr
0.2 & 0.8}, \qquad \bolds{\pi} = \pmatrix{
0.5
\cr
0.5},
\]
respectively. The emission densities $f_1$ and $f_2$ are skew-normal
densities [\citet{Azzalini1985}],
\[
f_{j}(y) = 2\phi(y-\mu_j) \Phi\bigl(
\alpha_j(y-\mu_j)\bigr),
\]
with $\mu_1=-2$, $\alpha_1=5$ and $\mu_2=-\mu_1$, $\alpha
_2=-\alpha_1$. The sign of the skewness parameters $\alpha_1,\alpha
_2$ implies that $f_1$ is skewed to the right while $f_2$ is skewed to
the left.

In each of $500$ Monte Carlo replications, we estimated the two
emission densities $f_1$ and $f_2$ using our orthogonal-series
estimator and constructed $95\%$ confidence intervals at the
percentiles of $f_1$ and $f_2$. We present results for $n=500$ (left
plot) and $n={}$5000 (right plot) graphically in Figure~\ref{figmchmmcomplete}. Results for additional sample sizes are
available in the supplementary material [\citet{BonhommeJochmansRobin2014b}].

\begin{figure}

\includegraphics{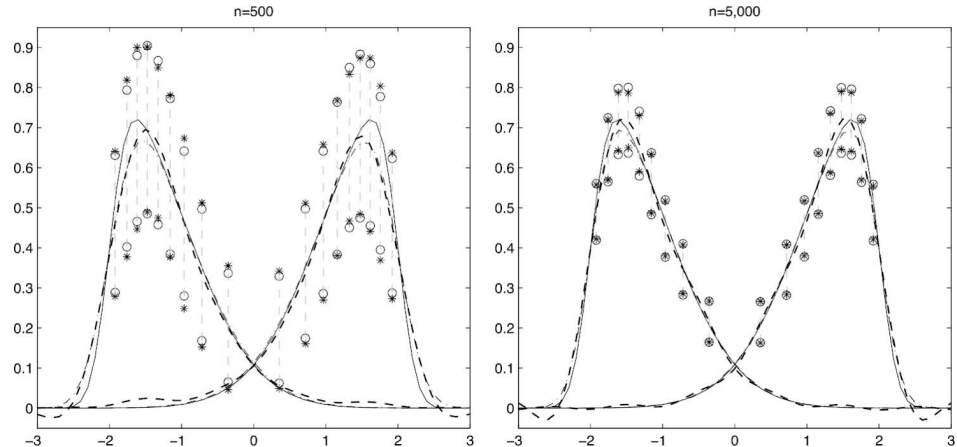}

\caption{Emission densities in the hidden Markov model.}
\label{figmchmmcomplete}
\end{figure}

Each plot in Figure~\ref{figmchmmcomplete} contains the true
functions $f_1$ and $f_2$ (solid lines), and the mean (across the Monte
Carlo replications) of our orthogonal-series estimator (dashed lines)
as well as of an infeasible kernel-density estimator (dashed--dotted
lines) computed from the subsample of observations that are in the
respective latent state (see the supplementary material for more
detail). The plots show that, even in small samples, our estimator
essentially coincides with the infeasible estimator, on average.

Figure~\ref{figmchmmcomplete} also contains average $95\%$
confidence intervals ($-\circ$), based on the pointwise distributional
result in Theorem~\ref{thmpointwise}, for the emission densities at
their respective percentiles. To assess the adequacy of our asymptotic
approximation, the plots in the figure also provide $95\%$ confidence
intervals at the percentiles constructed using the empirical standard
deviation of the point estimates across the Monte Carlo replications
($-\ast$). Figure~\ref{figmchmmcomplete} shows that our estimated
standard error captures well the small-sample variability of the
orthogonal-series estimator.


\section*{Acknowledgements}
We thank the Editor (Runze Li), an Associate Editor, three referees,
Xiaohong Chen, Ignat Domanov, Marc Henry and Nick Vannieuwenhoven for
comments. We are grateful to Laurent Albera and Xavier Luciani for
sharing the code for their diagonalization algorithm in \citeauthor{LucianiAlbera2010}
(\citeyear{LucianiAlbera2010,LucianiAlbera2014}) with us. Early versions of this
paper circulated as ``Nonparametric spectral-based estimation of latent
structures.''

\begin{supplement}[id=suppA]
\stitle{Supplement to ``Estimating multivariate latent-structure models''}
\slink[doi]{10.1214/15-AOS1376SUPP} 
\sdatatype{.pdf}
\sfilename{aos1376\_supp.pdf}
\sdescription{The supplement to this paper [\citet{BonhommeJochmansRobin2014b}]
contains additional details and discussion, omitted proofs and
additional simulation results.}
\end{supplement}


\printaddresses
\end{document}